\documentclass[12pt]{amsart}

\usepackage{amsmath}
\usepackage{amscd,amssymb}
\usepackage{epsfig}

\usepackage{graphicx}
\usepackage{color}
\usepackage{hyperref}
\usepackage{booktabs}
\usepackage{amsfonts}
\usepackage{mathrsfs}

\renewcommand{\proof}{\par\noindent{\it Proof.\ \ }}
\def\qed{\ifmmode\square\else\nolinebreak\hfill
$\Box$\fi\par\vskip12pt}

\newtheorem{thm}{Theorem}[section]
\newtheorem{lemma}[thm]{Lemma}
\newtheorem{corollary}[thm]{Corollary}

\numberwithin{equation}{section}
\numberwithin{thm}{section}

\theoremstyle{definition}

\newtheorem{remark}[thm]{Remark}

\newcommand{\bF}{\mathbb F}

\newcommand{\bZ}{\mathbb Z}

\newcommand{\cS}{\mathcal S}

\definecolor{Purple}{rgb}{0.5,0,0.5}

\def\Aut{{\rm Aut}}

     \def\rad{{\rm Rad}}
\def\Ker{{\rm Ker}}  \def\Tr{{\rm Tr}}

\newcommand{\<}{\langle}
\renewcommand{\>}{\rangle}

\begin{document}
\pagestyle{plain}
\begin{titlepage}

\title{Partial difference sets and amorphic Cayley schemes in non-abelian $2$-groups}
\begin{center}
\author{Tao Feng, Zhiwen He$^*$ and Yu Qing Chen
}\end{center}
\address{School of Mathematical Sciences, Zhejiang University, Hangzhou 310027, China}
\email{tfeng@zju.edu.cn}

\address{School of Mathematical Sciences, Zhejiang University, Hangzhou 310027,  China}
\email{zhiwen$\_$he@zju.edu.cn}

\address{Department of Mathematics and Statistics, Wright State University, Dayton, OH 45435-0001}

\email{yuqing.chen@wright.edu}

\begin{abstract} In this paper, we consider regular automorphism groups of graphs in the RT$2$ family and the Davis-Xiang family and amorphic abelian Cayley schemes from these graphs. We derive general results on the existence of non-abelian regular automorphism groups from abelian regular automorphism groups and apply them to the RT$2$ family and Davis-Xiang family and their amorphic abelian Cayley schemes to produce amorphic non-abelian Cayley schemes.
\end{abstract}

\keywords{amorphic Cayley scheme; Latin square type; negative Latin square type; partial difference set; regular automorphism group; strongly regular graph.\\
{\bf  Mathematics Subject Classification (2010) 05E30 05E18 05C25}\\
$^*$Correspondence author}

\maketitle

%%%%%%%%%%%%%%%%%%%%%%%%%%%%%%%%%%%

\section{Introduction}
 A strongly regular graph $\Gamma$ with parameters $(v,k,\lambda,\mu)$ is a regular graph of valency $k$, not complete or null, such that any two distinct vertices $x$ and $y$ have exactly $\lambda$ or $\mu$ common neighbors according as $x$ and $y$ are adjacent or not respectively. In the case $\Gamma$ has an automorphism group $G$ acting regularly on the vertex set, we fix any vertex $a$, and define $D=\{g\in G:a^g\text{ is adjacent to } a\}\subseteq G$. Then $\Gamma$ is isomorphic to the Cayley graph $\textup{Cay}(G,D)$, whose vertices are the elements of $G$, and two vertices $g_1$ and $g_2$ are adjacent if and only if $g_2g_1^{-1}\in D$. The fact that $\Gamma$ is loopless and undirected strongly regular is equivalent to $1\notin D$, $D=D^{(-1)}$, and $|D\cap gD|=\lambda$ or $\mu$ according as $g$ is a non-identity element in $D$ or $G\setminus D$ respectively, where $D^{(-1)}=\{g^{-1}:\,g\in D\}$ and $gD=\{gh:\,h\in D\}$. We call such a set $D$ with these properties a (non-trivial regular) partial difference set (PDS) with parameters $(v,k,\lambda,\mu)$ in the group $G$. We refer the reader to \cite{32} for more background on strongly regular graphs.

Partial difference sets have close connections with other branches of combinatorics as well as coding theory and finite geometry. The reader can find many partial difference sets constructed from projective two-weight codes and projective two-intersecting sets in \cite{27}. In \cite{31}, Ma provided a thorough and comprehensive survey on the development of the subject up to 1994. Since then many new construction methods and techniques have been discovered, polished, and perfected, including the use of  finite fields and finite local rings to construct partial difference sets in finite abelian $p$-groups (see \cite{3,12,7,11,4,22,9,10,20,2}). New projective two-weight codes in \cite{25,27,26} also yield new partial difference sets in elementary abelian $p$-groups.  We now try to shift our attention to partial difference sets in non-abelian groups. The only known constructions of such partial difference sets, to our best knowledge, are those in \cite{Ghinelli} by Ghinelli,  \cite{23} by Swartz and \cite{smith} by Smith.

In this paper, we consider partial difference set with parameter $(n^2,r(n-\epsilon),\epsilon n+r^2-3\epsilon r, r^2-\epsilon r)$ for $\epsilon=\pm 1$.  When $\epsilon=1$ (resp. $\epsilon=-1$), the partial difference set with such a parameter set is called a Latin square  type (resp. negative Latin square type) partial difference set, and we write LS (resp. NLS) type for short. There are many constructions of partial difference sets with both parameter sets (see \cite{28,21,33,24,6,16,15,8,17,1,5,34}). The known partial difference sets with NLS parameters are relatively rare. Almost all the known groups that contain LS or NLS type partial difference sets are abelian $p$-groups and many of these partial difference sets are obtained from quadratic forms \cite{24}, bent functions \cite{28,21,33,34} and various other combinatorial objects. Most notably, Davis and Xiang \cite{11} constructed the first known family of NLS type partial difference sets in non-elementary abelian 2-groups of exponent $4$ by using quadrics.  For product constructions, the reader may consult Polhill et al. \cite{13,19}.

A symmetric association scheme $\cS=(X; \Gamma_1, \Gamma_2, \cdots, \Gamma_n)$ on a finite set $X$ is a partition of the edge set of the complete graph on $X$ into spanning subgraphs $\Gamma_1$, $\Gamma_2$, $\cdots$, $\Gamma_n$ on $X$ such that their adjacency matrices $A_1$, $A_2$, $\cdots$, $A_n$ satisfy $A_iA_j=p_{i,j}^0I_X+p_{i,j}^1A_1+p_{i,j}^2A_2+\cdots+p_{i,j}^nA_n$ for some non-negative integers $p_{i,j}^0$, $p_{i,j}^1$, $\cdots$, $p_{i,j}^n$ for all $1\leqslant i,j\leqslant n$, where $I_X$ is the identity matrix. Such an association scheme is said to be amorphic if for any partition $\Gamma'_1$, $\Gamma'_2$, $\cdots$, $\Gamma'_{m}$ of the edge set of the complete graph on $X$ satisfying $\Gamma_i\subseteq\Gamma'_j$ for some $1\leqslant j\leqslant m$ for all $1\leqslant i\leqslant n$, $\cS'=(X; \Gamma'_1, \Gamma'_2, \cdots, \Gamma'_{m})$ is a symmetric association scheme on $X$. Combining the results of Ivanov \cite{Ivanov} and  of Ito, Munemasa and Yamada \cite{Ito}, we know that the symmetric association scheme $\cS$ with $n>2$ is amorphic if and only if the graphs $\Gamma_1$, $\Gamma_2$, $\cdots$, $\Gamma_n$ are all strongly regular graphs of LS type or all strongly regular graphs of NLS type. Furthermore, van Dam \cite{Dam} proved that if a partition of the edge set of the complete graph on $X$ into spanning subgraphs $\Gamma_1$, $\Gamma_2$, $\cdots$, $\Gamma_n$ such that they are all strongly regular graphs of LS type or all strongly regular graphs of NLS type, then $\cS=(X; \Gamma_1, \Gamma_2, \cdots, \Gamma_n)$ is a symmetric association scheme, and hence an amorphic association scheme. A permutation of $X$ is said to be an automorphism of $\cS$ if it is an automorphism of all subgraphs $\Gamma_1$, $\Gamma_2$, $\cdots$, $\Gamma_n$. The set of all automorphisms of $\cS$ forms a group, called the full automorphism group of $\cS$ and denoted by $\Aut(\cS)$. Any subgroup of $G$ of $\Aut(\cS)$ is called an automorphism group of $\cS$. We say that the association scheme $\cS=(X; \Gamma_1, \Gamma_2, \cdots, \Gamma_n)$ is a Cayley association scheme if $\cS$ has an automorphism group $G$ that acts regularly on the set $X$. When $G$ is a regular automorphism group of $\cS=(X; \Gamma_1, \Gamma_2, \cdots, \Gamma_n)$, there is a partition $D_1$, $D_2$, $\cdots$, $D_n$ of $G\setminus\{1\}$ such that $D_i=D_i^{(-1)}$ and $\Gamma_i={\rm Cay}(G,D_i)$ for all $1\leqslant i\leqslant n$, and such a Cayley association scheme is written $\cS=(G;D_1,D_2,\cdots,D_n)$. The Cayley association scheme $\cS=(G;D_1,D_2,\cdots,D_n)$ with $n>2$ is amorphic if and only if $D_1,D_2,\cdots,D_n$ are all partial difference sets of LS type in $G$ or all partial difference sets of NLS type in $G$. By a theorem of van Dam \cite{Dam}, any partition $D_1,D_2,\cdots,D_n$ of $G\setminus\{1\}$ into partial difference sets of LS type in $G$ or into partial difference sets of NLS type in $G$ yields a Cayley association scheme $\cS=(G;D_1,D_2,\cdots,D_n)$ and hence an amorphic Cayley association scheme. 

One of the central problems in the study of partial difference sets is that for a given parameter set, determine all groups of the appropriate order that contain a partial difference set with the given parameters, as noted in \cite{11}. In this paper, we consider regular subgroups of the graphs and Cayley association schemes arising from the RT2 family and the Davis-Xiang family, which will be defined in the next section. Such regular subgroups can be partitioned into partial difference sets with the same parameters as those of the families.  The regular groups we construct are non-abelian $2$-groups of nilpotency class $2$, $3$, $4$ or $6$ and of exponent $4$, $8$ or $16$.

This paper is organized as follows. In Section 2, we introduce a general framework to describe the RT2 family and the Davis-Xiang family and their amorphic association schemes, and describe our general strategy to study their regular subgroups. In Section 3, we first derive some general results on the existence of regular groups of partial difference set constructed from quadratic forms, and then apply them to the RT2 family and the Davis-Xiang family and amorphic association schemes. We also compute some invariants of those regular automorphism groups. We conclude the paper with Section 4.

\end{titlepage}
\section{Preliminaries}

\subsection{The RT2 and the Davis-Xiang graphs and their amorphic association schemes}

Let $q$ be a power of a prime and $\bF_q$ be the finite field with $q$ elements. Let $V$ be a vector space of dimension $n$ over $\bF_q$. A function $Q\colon V\longrightarrow\bF_q$ is a quadratic form if it satisfies:
\begin{itemize}
\item[(i)] $Q(sv)=s^2Q(v)$ for all $s\in \bF_q$ and $v\in V$;
\item[(ii)] the function $B(v_1,v_2):=Q(v_1+v_2)-Q(v_1)-Q(v_2)$ is $\bF_q$-bilinear.
\end{itemize}
The radical of a quadratic form $Q\colon V\longrightarrow\bF_q$ is
\[
\rad(Q):=\{v\in V:Q(v)=0\text{ and }B(v,x)=0\textup{ for all }x\in V\}.
\]
The quadratic form $Q$ is nonsingular if it has a trivial radical. If this is the case, then there is a basis $e_1,\cdots,e_n$ of $V$ such that for each element $x=\sum_{i=1}^nx_ie_i$ in $V$, exactly one of the following occurs:
\begin{enumerate}
\item $Q(x)=x_1x_2+\cdots+x_{n-1}x_n$, $n$ even (hyperbolic);
\item $Q(x)=x_1x_2+\cdots+x_{n-3}x_{n-2}+x_{n-1}^2+ax_{n-1}x_n+bx_n^2$  for some $a$ and $b$ in $\bF_q$ with $X^2+aX+b$ irreducible over $\bF_q$, $n$ even (elliptic);
\item $Q(x)=x_1x_2+\cdots+x_{n-2}x_{n-1}+x_n^2$, $n$ odd (parabolic).
\end{enumerate}
When $n$ is even and $Q$ is non-singular, we define $s(Q)=1$ when $Q$ is hyperbolic and $s(Q)=-1$ when $Q$ is elliptic. If $Q_1\colon V_1\longrightarrow\bF_q$ and $Q_2\colon V_2\longrightarrow\bF_q$ are two non-singular quadratic forms of even dimension, then $(Q_1\oplus Q_2)\colon V_1\oplus V_2\longrightarrow\bF_q$, given by $(Q_1\oplus Q_2)(v_1,v_2)=Q_1(v_1)+Q_2(v_2)$ for all vectors $v_1\in V_1$ and $v_2\in V_2$, is also a non-singular quadratic form of even dimension and $s(Q_1\oplus Q_2)=s(Q_1)s(Q_2)$. A bijective $\bF_q$-linear transformation $g:\,V\rightarrow V$ is an isometry of $Q$ if  $Q(g(x))=Q(x)$ for each $x\in V$. It is called a generalized isometry if $Q(g(x))=Q(x)^\sigma$ for each $x\in V$, where $\sigma$ is an element of the absolute Galois group ${\rm Gal}(\bF_q)$.

There is a classical construction of strongly regular graphs from the zeros of a non-singular quadratic form, which we call \textit{RT2 graphs} following \cite{27}.
{\thm\label{QD}
Let $V$ be a $2l$-dimensional vector space over $\bF_q$, equipped with a non-singular quadratic form $Q$. The set $Q^{-1}(0)\setminus\{0\}=\{x\in V:\,x\ne 0\text{ and }Q(x)=0\}$ is a partial difference set with parameters
\begin{equation*}
(q^{2l},(q^{l-1}+s(Q))(q^l-s(Q)),q^{2l-2}+s(Q)q^{l-1}(q-1)-2,q^{2l-2}+s(Q)q^{l-1}).
\end{equation*}
When $Q$ is elliptic, the partial difference set $Q^{-1}(0)\setminus\{0\}$ is of  NLS type, and when $Q$ is hyperbolic, the partial difference set $Q^{-1}(0)\setminus\{0\}$ is of LS type.}

In \cite{11}, Davis and Xiang discovered a very intriguing way of twisting the above construction to obtain partial difference sets in non-elementary abelian $2$-groups. We call the strongly  regular graphs arising from their construction \textit{Davis-Xiang graphs}. We now describe a unified construction of both RT2 graphs and Davis-Xiang graphs and their amorphic Cayley association schemes in certain $2$-groups.

Let $\bF_4=\{0,1,\omega,\omega+1\}$ be the finite field of $4$ elements, where $\omega$ is a primitive element of $\bF_4$. It is interesting to note that for any $\alpha$ and $\beta$ in $\bF_4$, the quadratic form $Q_{\alpha,\beta}\colon\bF_4\times\bF_4\longrightarrow\bF_4$ given by $Q_{\alpha,\beta}(x,y)=\alpha x^2+xy+\beta y^2$ has the property that $Q_{\alpha,\beta}^{-1}(x)$ is a  $(16,3,2,0)$-partial difference set for all $x\ne0$ in $\bF_4$ when $s(Q_{\alpha,\beta})=1$ and a $(16,5,0,2)$-partial difference set for all $x\ne0$ in $\bF_4$ when $s(Q_{\alpha,\beta})=-1$.
Let $\epsilon$ be an element in $\bF_2$. We define an abelian group $G_\epsilon$ whose underlying set is $\bF_4\times\bF_4$ and whose addition is $(x,y)+(x',y')=(x+x',y+y'+\epsilon(xx')^2)$ for all $(x,y)$ and $(x',y')\in\bF_4\times\bF_4$. When $\epsilon=0$, the group $G_\epsilon\cong\bZ_2^4$, and when $\epsilon=1$, the group $G_\epsilon\cong\bZ_4^2$. Let $\alpha$ be an element in $\bF_4$. We define a quadratic form $Q_\alpha\colon \bF_4\times\bF_4\longrightarrow\bF_4$ by $Q_\alpha(x,y)=\alpha x^2+xy+y^2$ for all $(x,y)\in\bF_4\times\bF_4$. The value $s(Q_\alpha)=(-1)^{{\rm Tr}(\alpha)}$, where ${\rm Tr}(\alpha)=\alpha+\alpha^2$ is the trace of $\alpha\in\bF_4$. The quadratic form $Q_\alpha$ can be regarded as a function $Q_\alpha\colon G_\epsilon\longrightarrow\bF_4$ for both $\epsilon=0$ and $\epsilon=1$. It is easy to check that $Q^{-1}_\alpha(x)$ is a $(16,4-s(Q_\alpha),1+s(Q_\alpha),1-s(Q_\alpha))$-partial difference set in both groups $G_0$ and $G_1$ for all $x\ne0$ in $\bF_4$; Consequently, $Q^{-1}_\alpha(0)\setminus\{0\}$ is a $(16,6,2,2)$-partial difference set, also known as a $(16,6,2)$-Hadamard difference set, in $G_0$ and $G_1$ when $s(Q_\alpha)=1$, and $Q^{-1}_\alpha(0)\setminus\{0\}=\emptyset$ for both $G_0$ and $G_1$ when $s(Q_\alpha)=-1$. The following Theorem is a consequence of Proposition 3.6 and Theorem 3.7 in \cite{21}.

{\thm\label{thm_DX}
Let $n\geqslant2$ be an integer. Given a vector $\mathbf{a}=(\alpha_1,\alpha_2,\cdots,\alpha_n)$ in $\bF_4^n$ and a vector $\mathbf{e}=(\epsilon_1,\epsilon_2,\cdots,\epsilon_n)$ in $\bF_2^n$, we define a quadratic form $Q_\mathbf{a}=Q_{\alpha_1}\oplus Q_{\alpha_2}\oplus\cdots\oplus Q_{\alpha_n}$ and an abelian group $G_\mathbf{e}=G_{\epsilon_1}\oplus G_{\epsilon_2}\oplus\cdots\oplus G_{\epsilon_n}$, and the quadratic form $Q_\mathbf{a}$ can be regarded as a function $Q_\mathbf{a}\colon G_\mathbf{e}\longrightarrow\bF_4$. The subset $Q_\mathbf{a}^{-1}(0)\setminus\{0\}$ is a $$(4^{2n},(4^{n-1}+s(Q_\mathbf{a}))(4^n-s(Q_\mathbf{a})),4^{2n-2}+3\cdot4^{n-1}s(Q_\mathbf{a})-2,4^{2n-2}+4^{n-1}s(Q_\mathbf{a}))$$
partial difference set in $G_\mathbf{e}$ for each $\mathbf{a}\in\bF_4^n$ and each $\mathbf{e}\in\bF_2^n$, and $Q_\mathbf{a}^{-1}(x)$ is a $$(4^{2n},4^{n-1}(4^n-s(Q_\mathbf{a})),4^{2n-2}+4^{n-1}s(Q_\mathbf{a}),4^{2n-2}-4^{n-1}s(Q_\mathbf{a}))$$
partial difference set in $G_\mathbf{e}$ for each $x\ne 0$ in $\bF_4$, each $\mathbf{a}\in\bF_4^n$ and each $\mathbf{e}\in\bF_2^n$. Clearly when $s(Q_\mathbf{a})=-1$, all partial difference sets $Q_\mathbf{a}^{-1}(0)\setminus\{0\}$, $Q_\mathbf{a}^{-1}(1)$, $Q_\mathbf{a}^{-1}(\omega)$ and $Q_\mathbf{a}^{-1}(\omega+1))$ are of NLS type, and when $s(Q_\mathbf{a})=1$, all partial difference sets $Q_\mathbf{a}^{-1}(0)\setminus\{0\}$, $Q_\mathbf{a}^{-1}(1)$, $Q_\mathbf{a}^{-1}(\omega)$ and $Q_\mathbf{a}^{-1}(\omega+1))$ are of LS type. Hence $$\cS^{(4)}_{\mathbf{e},\mathbf{a}}=(G_\mathbf{e}; Q_\mathbf{a}^{-1}(0)\setminus\{0\}, Q_\mathbf{a}^{-1}(1), Q_\mathbf{a}^{-1}(\omega), Q_\mathbf{a}^{-1}(\omega+1))$$ and $$\cS^{(3)}_{\mathbf{e},\mathbf{a}}=(G_\mathbf{e}; Q_\mathbf{a}^{-1}(0)\setminus\{0\}, Q_\mathbf{a}^{-1}(1), Q_\mathbf{a}^{-1}(\omega)\cup Q_\mathbf{a}^{-1}(\omega+1))$$are amorphic abelian Cayley schemes.}
{\remark When the vector $\mathbf{e}=\mathbf{0}$, the partial difference set $Q_\mathbf{a}^{-1}(0)\setminus\{0\}$ in $G_\mathbf{e}$ gives rise to RT2 graphs with $q=4$. When $\mathbf{e}\ne\mathbf{0}$, the graph generated by  partial difference set $Q_\mathbf{a}^{-1}(0)\setminus\{0\}$ in $G_\mathbf{e}$ is the Davis-Xiang graph.}
{\remark It is first proved by Davis and Xiang in \cite{11} that $Q_\mathbf{a}^{-1}(0)\setminus\{0\}$ is a partial difference set in $G_\mathbf{e}$. Using product method, Polhill \cite{20} showed that $Q_\mathbf{a}^{-1}(x)$ is a partial difference set in $G_\mathbf{e}$ for each $x\ne 0$ in $\bF_4$.}
{\remark The union of $Q_\mathbf{a}^{-1}(0)$ and $Q_\mathbf{a}^{-1}(1)$ and the union of $Q_\mathbf{a}^{-1}(\omega)$ and $Q_\mathbf{a}^{-1}(\omega+1)$ are Hadamard difference sets in $G_\mathbf{e}$.}

It is asked in \cite{11} that ``... one of the central problems in the study of partial difference sets is that for a given parameter set, which groups of the appropriate order contain a partial difference set with that parameter". Our construction of new regular automorphism groups is an effort to provide some answers to the question.

\subsection{Regular groups of strongly regular graphs and amorphic association schemes}

In this subsection, we give a brief description of our strategy to obtain regular automorphism groups of the Cayley schemes  $\cS^{(4)}_{\mathbf{e},\mathbf{a}}$ and $\cS^{(3)}_{\mathbf{e},\mathbf{a}}$ as defined in Theorem \ref{thm_DX} by constructing regular automorphism groups of the Cayley graphs from the partial difference sets of the schemes.  For each $g\in G_\mathbf{e}$, the map
\begin{align*}
&R(g)\colon G_\mathbf{e}\longrightarrow G_\mathbf{e}\\
&R(g)(x)=x+g\text{ for all }x\in G_\mathbf{e}
\end{align*}
preserves adjacency of strongly regular graphs in $\cS^{(4)}_{\mathbf{e},\mathbf{a}}$ and $\cS^{(3)}_{\mathbf{e},\mathbf{a}}$. Clearly $R(g_1)R(g_2)=R(g_1+g_2)$ for all $g_1,g_2\in G_\mathbf{e}$. For a subgroup $K$ of $G_\mathbf{e}$, we write $R(K)=\{R(g):g\in K\}$, which is a subgroup of $R(G)$ and is isomorphic to $K$.
{\lemma\label{lem_2ndAut} Let $\tau\in\Aut(G_\mathbf{e})$ be an automorphism of $G_\mathbf{e}$.
\begin{itemize}
\item[(i)] If $\tau$ is an isometry of $Q_\mathbf{a}$, i.e.,  $Q_\mathbf{a}(\tau(g))=Q_\mathbf{a}(g)$ for all $g\in G_\mathbf{e}$, then $\tau\in\Aut(\cS^{(4)}_{\mathbf{e},\mathbf{a}})$.
\item[(ii)] If $\tau$ is an generalized isometry of $Q_\mathbf{a}$, i.e.,  $Q_\mathbf{a}(\tau(g))=Q_\mathbf{a}(g)^2$ for all $g\in G_\mathbf{e}$, then $\tau\in\Aut(\cS^{(3)}_{\mathbf{e},\mathbf{a}})$.
\end{itemize}}
\proof Given a $\tau\in\Aut(G_\mathbf{e})$, the automorphism $\tau\in\Aut(\cS^{(3)}_{\mathbf{e},\mathbf{a}})$ or $\tau\in\Aut(\cS^{(4)}_{\mathbf{e},\mathbf{a}})$ if and only if $\tau(D)=D$ for every partial defference set $D$ from the amorphic abelian Cayley schemes $\cS^{(3)}_{\mathbf{e},\mathbf{a}}$ or $\cS^{(4)}_{\mathbf{e},\mathbf{a}}$. Hence if $Q_\mathbf{a}(\tau(g))=Q_\mathbf{a}(g)$ for all $g\in G_\mathbf{e}$, then $\tau\in\Aut(\cS^{(4)}_{\mathbf{e},\mathbf{a}})$, and if $Q_\mathbf{a}(\tau(g))=Q_\mathbf{a}(g)^2$ for all $g\in G_\mathbf{e}$, then $\tau\in\Aut(\cS^{(3)}_{\mathbf{e},\mathbf{a}})$.\qed

Lemma \ref{lem_2ndAut} ensures that $\<R(G_\mathbf{e}),\tau\>\cong G_\mathbf{e}\rtimes\<\tau\>$ is an automorphism group of $\cS^{(3)}_{\mathbf{e},\mathbf{a}}$ and $\cS^{(4)}_{\mathbf{e},\mathbf{a}}$. We will use $\tau$ and a subgroup of $G_\mathbf{e}$ to construct non-abelian regular automorphism groups of $\cS^{(3)}_{\mathbf{e},\mathbf{a}}$ and $\cS^{(4)}_{\mathbf{e},\mathbf{a}}$, which turn amorphic abelian Cayley schemes into amorphic non-abelian Cayley schemes.

{\thm\label{thm_GKt}
Suppose that  there is
\begin{enumerate}
\item[(a)] a generalized isometry $\tau\in\Aut(G_\mathbf{e})$ of $Q_\mathbf{a}$ of order $e>1$,
\item[(b)] a $\tau$-invariant subgroup $K$ of $G_\mathbf{e}$ of index $e$, and
\item[(c)] an element $h\in G_\mathbf{e}$ such that $h_e\in K$ and $h_1,\cdots,h_{e-1}\not\in K$, where
    \[
    h_i:=h+\tau(h)+\cdots+\tau^{i-1}(h),\quad 1\leqslant i\leqslant e.
    \]
\end{enumerate}
Then the group
\begin{equation}\label{eqn_GKt}
G_{K,\tau,h}:=\langle R(K),\,R(h)\tau\rangle
\end{equation}
is a regular automorphism group of the amorphic abelian Cayley scheme $\cS^{(3)}_{\mathbf{e},\mathbf{a}}$ or $\cS^{(4)}_{\mathbf{e},\mathbf{a}}$. Moreover, if we denote $G^{(0)}:=G_{K,\tau,h}$ and $G^{(k)}:=[G^{(k-1)}:G^{(0)}]$ for all $k\geqslant 1$, $Z(G^{(0)})$ the center of $G^{(0)}$ and $\Phi(G^{(0)})$ the Frattini subgroup of $G^{(0)}$, then
\begin{equation}\label{eqn_GKDer}
 G^{(1)}=\<R(x-\tau(x)):~x\in K\>\text{ and }G^{(k)}=\<R(x-\tau(x)):~x\in G^{(k-1)}\>\text{ for all }k\geqslant 2,
\end{equation}
and $G^{(0)}$ is non-abelian if and only if $\tau$ does not fix all elements of $K$, i.e. the order $t=o(\tau|_K)$ is greater than $1$. Let $m$ be the order of $R(h_t)\tau^t$, The center
\begin{align}
\label{eqn_center} Z(G^{(0)})\cong\begin{cases}[\left\langle R(x):\,x\in K\text{ with }\tau(x)=x\right\rangle/\left\langle R(h_e)\right\rangle]\times\bZ_m &\text{ if }t<e,\\
                                                \left\langle R(x):\,x\in K\text{ with }\tau(x)=x\right\rangle &\text{ if }t=e,\end{cases}
\end{align}
and the Frattini subgroup
\begin{equation}\label{eqn_Frat}
\Phi(G^{(0)})=\left\langle R(2x), R(x-\tau(x)), R(h_{2})\tau^{2}:\,x\in K\right\rangle=\left\langle R(2x), R(x+\tau(x)), R(h_{2})\tau^{2}:\,x\in K\right\rangle.
\end{equation}}
\proof It is routine to check that $(R(h)\tau)^{-1}=R(\tau^{-1}(h^{-1}))\tau^{-1}$ and
\[
 R(h)\tau R(x)\left(R(h)\tau\right)^{-1}=R(\tau(x))\quad\textup{ for } x\in K.
\]
It follows that $R(K)$ is a normal subgroup of $G_{K,\tau,h}$. Moreover,
we have $\left(R(h)\tau\right)^{i}=h_{i}\tau^i$ for $0\leqslant i\leqslant e$, where $h_0=0$ and the other $h_i$'s are defined in item (c). Since $h_e\in K$, we have
\begin{equation}\label{eqn_GKdec}
G_{K,\tau,h}=\cup_{i=1}^{e}R(K+h_i)\tau^i.
\end{equation}
In particular, $G_{K,\tau,h}$ has order $|G_\mathbf{e}|$.

For $0\leqslant i<j\leqslant e-1$, we have $h_jh_i^{-1}=\tau^i(h_{j-i+1})$, which is not in $K$ since $K$ is $\tau$-invariant and $h_{j-i}\not\in K$ by item (c). Together with $h_e\in K$, it follows that  $\cup_{i=1}^{e}(K+ h_i)=G_\mathbf{e}$. In order to show that $G_{K,\tau,h}$ is a regular subgroup, it suffices to prove that the orbit of $0\in G_\mathbf{e}$ under its action is  $G_\mathbf{e}$. Since $R(k)R(h_i)\tau^i(0)=h_i+k$, this is equivalent to $\cup_{i=1}^{e}(K+h_i)=G_\mathbf{e}$ which has just been established. This proves the first part. It is straightforward to compute that
\begin{align*}
G^{(1)}&=\< R(x)R(h_t)\tau^{t}R(x)^{-1}(R(h_t)\tau^{t})^{-1}:\, x\in K,\,1\leqslant t\leqslant e-1\>\\
             &=\<R(x-\tau^t(x)):\, x\in K,~1\leqslant t\leqslant e-1\>.
\end{align*}
Since $$x-\tau^t(x)=(x-\tau(x))+(\tau(x)-\tau(\tau(x)))+(\tau^2(x)-\tau(\tau^2(x)))+\cdots+(\tau^{t-1}(x)-\tau(\tau^{t-1}(x)))$$ and $K$ is $\tau$-invariant,
we get $$G^{(1)}=\<R(x-\tau^t(x)):\, x\in K,~1\leqslant t\leqslant e-1\>=\<R(x-\tau(x)):\, x\in K\>.$$ Similarly $G^{(k)}=\<R(x-\tau(x)):~x\in G^{(k-1)}\>$ for all $k\geqslant 2$.

It is  obvious that the center
\begin{align*}
Z(G_{K,\tau,h})&=\left\langle R(x),\,R(h_t)\tau^t:\,x\in K\textup{ with }\tau(x)=x,~t=o(\tau|_K)\right\rangle\\
               &\cong\begin{cases}[\left\langle R(x):\,x\in K\text{ with }\tau(x)=x\right\rangle/\left\langle R(h_e)\right\rangle]\times\bZ_m &\text{ if }t<e,\\
                                                \left\langle R(x):\,x\in K\text{ with }\tau(x)=x\right\rangle &\text{ if }t=e,\end{cases}
\end{align*}
as $\left\langle R(x): x\in K\text{ with }\tau(x)=x\right\rangle\cap\left\langle R(h_t)\tau^t\right\rangle=\left\langle R(h_e)\right\rangle$ and the order $o(R(h_t)\tau^t)=m$. Since $G_\mathbf{e}$ is a $2$-group, the Frattini subgroup $\Phi(G_{K,\tau,t})$ of $G_{K,\tau,t}$ is generated by $g^2$, $g\in G_{K,\tau,t}$. Since
\begin{align*}
R(x+h_i)\tau^iR(x+h_i)\tau^i&=R(x+h_i+\tau^i(x+h_i))\tau^{2i}=R(x+\tau^i(x)+h_{2i})\tau^{2i},
\end{align*}
we have
\begin{align*}
\Phi(G_{K,\tau,h})&=\left\langle R(x+\tau^i(x)+h_{2i})\tau^{2i}:\,x\in K,\,i\geqslant 0\right\rangle\\
                  &=\left\langle R(x+\tau^i(x)),~R(h_{2i})\tau^{2i}:\,x\in K,\,i\geqslant 0\right\rangle\\
                  &=\left\langle R(2x),~R(x+\tau^i(x)),~R(h_{2i})\tau^{2i}:\,x\in K,\,i\geqslant 1\right\rangle\\
                  &=\left\langle R(2x),~R(x-\tau^i(x)),~R((h_{2})\tau^{2})^i:\,x\in K,\,i\geqslant 1\right\rangle\\
                  &=\left\langle R(2x),~R(x-\tau(x)),~R(h_{2})\tau^{2}:\,x\in K\right\rangle\\
                  &=\left\langle R(2x),~R(x+\tau(x)),~R(h_{2})\tau^{2}:\,x\in K\right\rangle.
\end{align*}\qed
We shall use Eqns. \eqref{eqn_GKDer}, \eqref{eqn_center} and \eqref{eqn_Frat} as invariants to distinguish isomorphism classes of various regular subgroups $G_{K,\tau,h}$'s.

\section{Regular subgroups of the RT2 graphs and the Davis-Xiang graphs and their amorphic abelian Cayley schemes}
Given an abelian group $G$, we write ${\rm End}(G)$ for all endomorphisms of $G$. For each $f\in{\rm End}(G)$, we write ${\rm Im}_G(f):=\{f(g):g\in G\}$ for the image of $G$ under endomorphism $f$, and write ${\rm Ker}_G(f):=\{g\in G:f(g)=0\}$ for the kernel of endomorphism $f$ in $G$. For each integer $n$, we write ${\rm Im}_G(n):=\{ng:g\in G\}$ and ${\rm Ker}_G(n):=\{g\in G:ng=0\}$. For instance, the Frattini subgroup of $G_\mathbf{e}$ is $\Phi(G_\mathbf{e})={\rm Im}_{G_\mathbf{e}}(2)$ while ${\rm Fix}(\tau)={\rm Ker}_{G_\mathbf{e}}(\tau-1)$. Throughout this section, we assume that $G_\mathbf{e}$ is the abelian $2$-group defined in Theorem~\ref{thm_DX} and $\tau\in\Aut(G_\mathbf{e})$ is a generalized isometry of  $Q_\mathbf{a}$ of order $e=2$ or $4$, where $Q_\mathbf{a}$ is the quadratic form defined in Theorem~\ref{thm_DX}.

The main objective of this section is to find subspaces $K$ of $G_\mathbf{e}$ and elements $h$ in $G_\mathbf{e}$ that satisfy the requirements in Theorem \ref{thm_GKt} for the generalized isometry $\tau$. The following simple observation is very useful.
{\lemma\label{invariant}
The regular subgroup $G_{K,\tau,h}$ in Theorem \ref{thm_GKt} is non-abelian when $[G_\mathbf{e}:\,{\rm Fix}(\tau)]>e$.}
\proof By the second part of Theorem \ref{thm_GKt}, $G_{K,\tau,h}$ is abelian if and only if $\tau$ fixes each element of $K$, i.e., $K\leqslant{\rm Fix}(\tau)$. This is impossible by comparing their sizes.\qed

\subsection{General results}
Let $G$ be an abelian group and $\phi$ be an automorphism of $G$. We first prove some general results concerning certain $\phi$-invariant subgroups of $G$.
{\lemma\label{lem_index2}
Let $K$ be a subgroup of index $2$ in $G$. The subgroup $K$ is $\phi$-invariant if and only if ${\rm Im}_{G}(1+\phi)\leqslant K$.}
\proof If $K$ is $\phi$-invariant, then $x\in K$ if and only if $\phi(x)\in K$. Hence for any $x\notin K$, the element $(1+\phi)(x)=x+\phi(x)\in K$ as $[G:K]=2$. It is also clear that $(1+\phi)(x)=x+\phi(x)\in K$ when $x\in K$ as $K$ is a $\phi$-invariant subgroup. Therefore ${\rm Im}_{G}(1+\phi)\leqslant K$. On the other hand, if ${\rm Im}_{G}(1+\phi)\leqslant K$, then for any $x\in K$, the element $(1+\phi)(x)=x+\phi(x)\in K$, and hence $\phi(x)\in K$ and $K$ is $\phi$-invariant. \qed
{\lemma\label{lem_index4}
Let $K$ be a subgroup of index $4$ in $G$. If the order of $\phi$ is not divisible by $3$ and $K$ is $\phi$-invariant, then $G$ contains a subgroup $H$ of index $2$ such that $H$ is $\phi$-invariant and $K\leqslant H$. As a consequence, ${\rm Im}_{G}(1+\phi+\phi^2+\phi^3)\leqslant K$.}
\proof Since $K$ is a $\phi$-invariant subgroup of $G$, there is an induced automorphism $\bar{\phi}$ of the quotient group $\bar{G}:=G/K$. Let $\pi:G\longrightarrow\bar{G}$ be the natural projection homomorphism. Since the order of $\phi$ is not divisible by $3$ and $|\bar{G}|=4$, there must be an element $x$ of order $2$ in $\bar{G}$ such that $\bar{\phi}(x)=x$. Let $H=\pi^{-1}(\{0,x\})$. Then $H$ is $\phi$-invariant, $[G:H]=2$ and $K\leqslant H$. By Lemma~\ref{lem_index2}, we have ${\rm Im}_{G}(1+\phi)\leqslant H$. Since $K$ is also $\phi^2$ invariant, Lemma~\ref{lem_index2} implies ${\rm Im}_{H}(1+\phi^2)\leqslant K$. Therefore
${\rm Im}_{G}(1+\phi+\phi^2+\phi^3)={\rm Im}_{G}((1+\phi^2)(1+\phi))\leqslant K$.\qed

We have a complete answer to the case $o(\tau)=2$.
{\thm\label{thm_Order2}
Let $\tau\in\Aut(G_\mathbf{e})$ be a generalized isometry or an isometry of order $2$ of $Q_\mathbf{a}$. There is a $\tau$-invariant subgroup $K$ in $G_\mathbf{e}$ of index $2$ such that $G_{K,\tau,h}$ as defined in Eqn. \eqref{eqn_GKt} is a regular subgroup of the amorphic abelian Cayley scheme $\cS^{(3)}_{\mathbf{e},\mathbf{a}}$ or $\cS^{(4)}_{\mathbf{e},\mathbf{a}}$  for any element $h\in G_\mathbf{e}\setminus K$. The regular subgroup $G_{K,\tau,h}$ is of nilpotency class $1$ if $\tau|_K=1_K$, or nilpotency class $2$ if $\tau|_K\ne1_K$ and $\tau|_{\Phi(K)}=1_{\Phi(K)}$, or nilpotency class $3$ if $\tau|_{\Phi(K)}\ne1_{\Phi(K)}$.}
\proof We know that the order of $\tau$ is $2$, therefore $|G_\mathbf{e}|\geqslant 4$. The number of maximal subgroups in $G_\mathbf{e}$ is $|G_\mathbf{e}/\Phi(G_\mathbf{e})|-1$, which is odd. Therefore there must be a $\tau$-invariant maximal subgroup $K$ of $G_\mathbf{e}$. Hence by Theorem~\ref{thm_GKt} and Lemma~\ref{invariant}, the group $G_{K,\tau,h}$ defined by Eqn. \eqref{eqn_GKt} is a regular automorphism group of the amorphic abelian Cayley scheme $\cS^{(3)}_{\mathbf{e},\mathbf{a}}$ or $\cS^{(4)}_{\mathbf{e},\mathbf{a}}$  for any element $h\in G_\mathbf{e}\setminus K$. If $c$ is the least positive integer such that $(1-\tau)^c|_K=0_K$, then by Theorem~\ref{thm_GKt}, the regular subgroup $G_{K,\tau,h}$ is of nilpotency class $c$. Since $(1-\tau)^2=2(1-\tau)$ and $(1-\tau)^3=0$, the group $G_{K,\tau,h}$ is of nilpotency class $1$ if $\tau|_K=1_K$, and nilpotency class $2$ if $\tau|_K\ne1_K$ and $\tau|_{\Phi(K)}=1_{\Phi(K)}$, and nilpotency class $3$ if $\tau|_{\Phi(K)}\ne1_{\Phi(K)}$.\qed

We next consider the case $o(\tau)=4$.

{\thm\label{thm_Order4}
Let $\tau\in\Aut(G_\mathbf{e})$ be a generalized isometry or an isometry of order $4$ of $Q_\mathbf{a}$. There is a $\tau$-invariant subgroup $K$ in $G_\mathbf{e}$ of index $4$ and an element $h\in G_\mathbf{e}$ such that $G_{K,\tau,h}$ as defined in Eqn. \eqref{eqn_GKt} is a regular automorphism group of the amorphic abelian Cayley scheme $\cS^{(3)}_{\mathbf{e},\mathbf{a}}$ or $\cS^{(4)}_{\mathbf{e},\mathbf{a}}$ if and only if there is a subgroup $H$ in $G_\mathbf{e}$ of index $2$ such that $${\rm Ker}_{G_\mathbf{e}}(1+\tau)+{\rm Im}_{G_\mathbf{e}}(1+\tau)\leqslant H\text{ and }{\rm Im}_{G_\mathbf{e}}(1+\tau)\not\leqslant\Phi(H)+{\rm Im}_{H}(1+\tau).$$ If such a regular subgroup $G_{K,\tau,h}$ exist, then its nilpotency class is no larger than $6$.}
\proof Suppose $K$ and $h$ as state in the theorem exist. By Lemma~\ref{lem_index4}, there is a $\tau$-invariant subgroup $H$ in $G_\mathbf{e}$ of index $2$ such that $K\leqslant H$. By Lemma~\ref{lem_index2}, we have ${\rm Im}_{G_\mathbf{e}}(1+\tau)\leqslant H$, ${\rm Im}_{H}(1+\tau)\leqslant K$ and $\Phi(H)\leqslant K$ as $[G_\mathbf{e}:H]=[H:K]=2$, which also implies $$\Phi(H)+{\rm Im}_{H}(1+\tau)\leqslant K.$$ Since $(1+\tau)(h)\notin K$ by (c) in Theorem~\ref{thm_GKt}, and $(1+\tau)(h)\in H$, we find that $${\rm Im}_{H}(1+\tau)<{\rm Im}_{G_\mathbf{e}}(1+\tau)\text{ and }{\rm Im}_{G_\mathbf{e}}(1+\tau)\not\leqslant \Phi(H)+{\rm Im}_{H}(1+\tau)\leqslant K.$$ We now show that ${\rm Ker}_{G_\mathbf{e}}(1+\tau)\leqslant H$
by contradiction. Suppose ${\rm Ker}_{G_\mathbf{e}}(1+\tau)\not\leqslant H$. Then ${\rm Ker}_H(1+\tau)={\rm Ker}_{G_\mathbf{e}}(1+\tau)\cap H$ is an index $2$ subgroup of ${\rm Ker}_{G_\mathbf{e}}(1+\tau)$ as $[G_\mathbf{e}:H]=2$. From $G_\mathbf{e}/{\rm Ker}_{G_\mathbf{e}}(1+\tau)\cong{\rm Im}_{G_\mathbf{e}}(1+\tau)$ and $H/{\rm Ker}_H(1+\tau)\cong{\rm Im}_H(1+\tau)$, we find that $|{\rm Im}_H(1+\tau)|=|{\rm Im}_{G_\mathbf{e}}(1+\tau)|$, which implies that ${\rm Im}_H(1+\tau)={\rm Im}_{G_\mathbf{e}}(1+\tau)$ contradicting ${\rm Im}_{H}(1+\tau)<{\rm Im}_{G_\mathbf{e}}(1+\tau)$. Hence ${\rm Ker}_{G_\mathbf{e}}(1+\tau)+{\rm Im}_{G_\mathbf{e}}(1+\tau)\leqslant H$.

Conversely, if $H$ is an index $2$ subgroup of $G_\mathbf{e}$ such that ${\rm Ker}_{G_\mathbf{e}}(1+\tau)+{\rm Im}_{G_\mathbf{e}}(1+\tau)\leqslant H$ and ${\rm Im}_{G_\mathbf{e}}(1+\tau)\not\leqslant\Phi(H)+{\rm Im}_{H}(1+\tau)$, by Lemma~\ref{lem_index2}, the subgroup $H$ is $\tau$-invariant and there is an element $h\in G_\mathbf{e}\setminus H$ such that $h+\tau(h)=(1+\tau)(h)\notin\Phi(H)+{\rm Im}_{H}(1+\tau)$. Hence $H$ has a index $2$ subgroup $K$ such that $h+\tau(h)\notin K$ and ${\rm Im}_{H}(1+\tau)\leqslant K$, and by Lemma~\ref{lem_index2}, $K$ is $\tau$-invariant in $H$. Since $H$ is $\tau$-invariant in $G_\mathbf{e}$, $K$ is also $\tau$-invariant in $G_\mathbf{e}$ and $[G_\mathbf{e}:K]=4$.
By Lemma~\ref{lem_index4}, $h+\tau(h)+\tau^2(h)+\tau^3(h)=(1+\tau+\tau^2+\tau^3)(h)\in{\rm Im}_{G_\mathbf{e}}(1+\tau+\tau^2+\tau^3)\leqslant K$. From $h\notin H$ and $h+\tau(h)\in H$, the element $h+\tau(h)+\tau^2(h)\notin H$ as $H$ is $\tau$-invariant in $G_\mathbf{e}$. Hence $K$ and $h$ satisfy (a), (b) and (c) in Theorem~\ref{thm_GKt}.

Since $(1-\tau)^6=0$ as $\tau^4=1$ and $4x=0$ for all $x\in G_\mathbf{e}$, by Theorem~\ref{thm_GKt}, the nilpotency class of $G_{K,\tau,h}$ is less than or equal to $6$.\qed

{\corollary Let $\tau\in\Aut(G_\mathbf{e})$ be a generalized isometry or an isometry of order $4$ of $Q_\mathbf{a}$ such that the intersection ${\rm Im}_{G_\mathbf{e}}(1+\tau)\cap\Phi(G_\mathbf{e})=\{0\}$.  There is a $\tau$-invariant subgroup $K$ in $G_\mathbf{e}$ of index $4$ and an element $h\in G_\mathbf{e}$ such that $G_{K,\tau,h}$ as defined in Eqn. \eqref{eqn_GKt} is a regular automorphism group of the amorphic abelian Cayley scheme $\cS^{(3)}_{\mathbf{e},\mathbf{a}}$ or $\cS^{(4)}_{\mathbf{e},\mathbf{a}}$ if and only if the intersection ${\rm Ker}_{G_\mathbf{e}}(1+\tau)\cap{\rm Im}_{G_\mathbf{e}}(1+\tau)\ne\{0\}$, or equivalently, ${\rm Ker}_{G_\mathbf{e}}(1+\tau)<{\rm Ker}_{G_\mathbf{e}}(1+\tau)^2$. The nilpotency class of $G_{K,\tau,h}$ is no larger than $6$.}
\proof If ${\rm Ker}_{G_\mathbf{e}}(1+\tau)\cap{\rm Im}_{G_\mathbf{e}}(1+\tau)\ne\{0\}$, then ${\rm Ker}_{G_\mathbf{e}}(1+\tau)+{\rm Im}_{G_\mathbf{e}}(1+\tau)\ne G_\mathbf{e}$ and there are subgroups $H$ of index $2$ in $G_\mathbf{e}$ such that ${\rm Ker}_{G_\mathbf{e}}(1+\tau)+{\rm Im}_{G_\mathbf{e}}(1+\tau)\leqslant H$. From the intersection ${\rm Im}_{G_\mathbf{e}}(1+\tau)\cap\Phi(G_\mathbf{e})=\{0\}$, we get ${\rm Im}_{G_\mathbf{e}}(1+\tau)\not\leqslant\Phi(H)+{\rm Im}_{H}(1+\tau)$ for every such index $2$ subgroup $H$ in $G_\mathbf{e}$. By Theorem~\ref{thm_Order4}, there is a $\tau$-invariant subgroup $K$ in $G_\mathbf{e}$ of index $4$ and an element $h\in G_\mathbf{e}$ such that $G_{K,\tau,h}$ as defined in Eqn. \eqref{eqn_GKt} is a regular subgroup of the amorphic abelian Cayley scheme $\cS^{(3)}_{\mathbf{e},\mathbf{a}}$ or $\cS^{(4)}_{\mathbf{e},\mathbf{a}}$. Conversely, if such a $\tau$-invariant subgroup $K$ in $G_\mathbf{e}$ of index $4$ and such an element $h\in G_\mathbf{e}$ exist, then Theorem~\ref{thm_Order4} implies that there is a subgroup $H$ in $G_\mathbf{e}$ of index $2$ such that ${\rm Ker}_{G_\mathbf{e}}(1+\tau)+{\rm Im}_{G_\mathbf{e}}(1+\tau)\leqslant H$, and hence ${\rm Ker}_{G_\mathbf{e}}(1+\tau)\cap{\rm Im}_{G_\mathbf{e}}(1+\tau)\ne\{0\}$.\qed

We prove one more result on the existence of $\tau$-invariant subspace $K$ of index $4$ in $G_\mathbf{e}$ stated in Theorem \ref{thm_Order2}.

{\thm Let $\tau\in\Aut(G_\mathbf{e})$ be a generalized isometry or an isometry of order $4$ of $Q_\mathbf{a}$. If the induced action of $\tau$ on $G_\mathbf{e}/(\Phi(G_\mathbf{e})+{\rm Im}_{G_\mathbf{e}}(1+\tau^2))$ is non-trivial, then there exists a pair $(K,\,h)$, and hence a regular automorphism group $G_{K,\tau,h}$ of the amorphic abelian Cayley scheme $\cS^{(3)}_{\mathbf{e},\mathbf{a}}$ or $\cS^{(4)}_{\mathbf{e},\mathbf{a}}$, that satisfies Theorem \ref{thm_GKt}. The nilpotency class of $G_{K,\tau,h}$ is no larger than $6$.}
\proof By Lemma~\ref{lem_index2},
\begin{align*} &\Phi(G_\mathbf{e})+{\rm Im}_{G_\mathbf{e}}(1+\tau^2)=\text{ the intersection of all $\tau^2$-invariant subgroups in $G_\mathbf{e}$ of index $2$, and }\\
&\Phi(G_\mathbf{e})+{\rm Im}_{G_\mathbf{e}}(1+\tau)=\text{ the intersection of all $\tau$-invariant subgroups in $G_\mathbf{e}$ of index $2$.}
\end{align*}
Hence $\Phi(G_\mathbf{e})+{\rm Im}_{G_\mathbf{e}}(1+\tau^2)\leqslant \Phi(G_\mathbf{e})+{\rm Im}_{G_\mathbf{e}}(1+\tau)$ as $\tau$-invariant subgroups are also $\tau^2$-invariant. Since the induced action of $\tau$ on $G_\mathbf{e}/(\Phi(G_\mathbf{e})+{\rm Im}_{G_\mathbf{e}}(1+\tau^2))$ is non-trivial, the subgroup ${\rm Im}_{G_\mathbf{e}}(1+\tau)\not\leqslant\Phi(G_\mathbf{e})+{\rm Im}_{G_\mathbf{e}}(1+\tau^2)$ and there is an element $h\in G_\mathbf{e}$ such that $h+\tau(h)\in{\rm Im}_{G_\mathbf{e}}(1+\tau)$ has a non-zero image in $G_\mathbf{e}/(\Phi(G_\mathbf{e})+{\rm Im}_{G_\mathbf{e}}(1+\tau^2))$. This implies that there is an index $2$ subgroup $H$ of $G_\mathbf{e}$ such that $h+\tau(h)\notin H$ as $h+\tau(h)\notin\Phi(G_\mathbf{e})$, and by Lemma~\ref{lem_index2}, $\tau(H)\ne H$ and $\tau^2(H)=H$. Let $K=H\cap\tau(H)$. Then $K$ is an index $4$ $\tau$-invariant subgroup of $G_\mathbf{e}$. From $h+\tau(h)\notin H$, it implies $h+\tau(h)\notin K$, and hence $h\notin K$ as $K$ is $\tau$-invariant. By Lemma~\ref{lem_index4}, $h+\tau(h)+\tau^2(h)+\tau^3(h)\in K$ and $h+\tau(h)+\tau^2(h)\notin K$ as $h\notin K$ and $K$ is $\tau$-invariant.\qed

The rest of this section is devoted to find non-abelian regular automorphism groups of RT2 graphs and Davis-Xiang graphs as well as amorphic abelian Cayley association schemes related to these graphs in Theorem \ref{thm_GKt}. The generalized isometry or isometry $\tau$ we shall use will have order $2$ or $4$, so the results in this subsection apply.

\subsection{Non-abelian regular automorphism groups of amorphic abelian Cayley schemes } For any two vectors $\mathbf{u}=(u_1,u_2,\cdots,u_n)$ and $\mathbf{v}=(v_1,v_2,\cdots,v_n)$ we define $$\mathbf{u}*\mathbf{v}:=(u_1v_1,u_2v_2,\cdots,u_nv_n)\text{ and }\Tr(\mathbf{u})=(\Tr(u_1),\Tr(u_2),\cdots,\Tr(u_n)).$$ We also denote by $\mathbf{1}:=(1,1,\cdots,1)$ and $\mathbf{0}:=(0,0,\cdots,0)$. We write $w(\mathbf{u})$ for the Hamming weight of the vector $\mathbf{u}$.

We start with the group $G_{\epsilon}$ with $\epsilon\in\bF_2$ and the quadratic form $Q_{\alpha}$ with $\alpha\in\bF_4$ defined before Theorem~\ref{thm_DX}. For each $\nu\in\bF_4$, we define two automorphisms $\rho_\nu\in\Aut(G_{\epsilon})$ and $\tau_\nu\in\Aut(G_{\epsilon})$ by $\rho_\nu(x,y)=(x^2,y^2+\nu x^2)$ and $\tau_\nu(x,y)=(x,y+\nu x)$ for all $(x,y)\in G_{\epsilon}$. The automorphism $\tau_\nu$ is an isometry of $Q_\alpha$ of order $1$ when $\nu=0$ and an isometry of $Q_\alpha$ of order $2$ when $\nu=1$. It is easy to check that $\rho_\nu^2=\tau_{\Tr(\nu)}$. The automorphism $\rho_\alpha$ is an generalized isometry of $Q_\alpha$ of order $2$ when $\Tr(\alpha)=0$ and a generalized isometry of $Q_\alpha$ of order $4$ when $\Tr(\alpha)=1$. It is easy to check that
\begin{align}    \notag&{\rm Im}_{G_\epsilon}(1+\tau_\nu)=\{(0,(\nu+\epsilon)x):x\in\bF_4\},\\
\label{o2}                           &{\rm Im}_{G_\epsilon}(1-\tau_\nu)=\{(0,\nu x):x\in\bF_4\},\\
                 \notag&\Ker_{G_\epsilon}(1-\tau_\nu)=\{((\nu+1)x,y)\in G_\epsilon: x,y\in\bF_4\},\end{align} and
\begin{align}
{\rm Im}_{G_\epsilon}(1+\rho_\alpha)\label{image_1+rho}&=\{(x,y)+\rho_\alpha(x,y): (x,y)\in G_\epsilon\}\\
                                  \notag  &=\{(x+x^2,y+y^2+\alpha x^2+\epsilon x^3): x,y\in\bF_4\}\\
                                 \notag   &=\begin{cases}\{(0,0),(0,1),(1,0),(1,1)\} & \text{ if } \alpha=0,\\ \{(0,0),(0,1),(1,\omega),(1,\omega+1)\} & \text{ if } \alpha=1,\\\{(0,y),(1,y): y\in\bF_4\} &\text{ if } \alpha=\omega\text { or }\omega+1.\end{cases}\\
{\rm Ker}_{G_\epsilon}(1+\rho_\alpha)\label{kernel_1+rho}&=\{(x,y)\in G_\epsilon : x+x^2=0 \text{ and } y+y^2+\alpha x^2+\epsilon x^3=0\text{ with }x,y\in\bF_4\}\\
                             \notag         &=\begin{cases}\{(0,0),(0,1),(1,0),(1,1)\}& \text{ if } \alpha+\epsilon=0,\\ \{(0,0),(0,1),(1,\omega),(1,\omega+1)\} &\text{ if } \alpha+\epsilon=1,\\ \{(0,0),(0,1)\} &\text{ if } \alpha=\omega\text { or }\omega+1.\end{cases}\\
{\rm Im}_{G_\epsilon}(1-\rho_\alpha)\label{image_1-rho}&=\{(x,y)-\rho_\alpha(x,y): (x,y)\in G_\epsilon\}\\
                                   \notag &=\{(x+x^2,y+y^2+(\alpha+\epsilon)x^2+\epsilon x^3): x,y\in\bF_4\}\\
                                   \notag &=\begin{cases}\{(0,0),(0,1),(1,0),(1,1)\} & \text{ if } \alpha+\epsilon=0,\\ \{(0,0),(0,1),(1,\omega),(1,\omega+1)\} & \text{ if } \alpha+\epsilon=1,\\\{(0,y),(1,y): y\in\bF_4\} &\text{ if } \alpha=\omega\text { or }\omega+1.\end{cases}\\
{\rm Ker}_{G_\epsilon}(1-\rho_\alpha)\label{kernel_1-rho}&=\{(x,y)\in G_\epsilon : x+x^2=0 \text{ and } y+y^2+\alpha x^2=0\text{ with }x,y\in\bF_4\}\\
                                    \notag &=\begin{cases}\{(0,0),(0,1),(1,0),(1,1)\}& \text{ if } \alpha=0,\\ \{(0,0),(0,1),(1,\omega),(1,\omega+1)\} &\text{ if } \alpha=1,\\ \{(0,0),(0,1)\} &\text{ if } \alpha=\omega\text { or }\omega+1.\end{cases}
\end{align}
We are now ready to present four families of amorphic non-abelian Cayley association schemes that are isomorphic to $\cS^{(3)}_{\mathbf{e},\mathbf{a}}$ or $\cS^{(4)}_{\mathbf{e},\mathbf{a}}$, and consequently non-abelian partial difference sets of LS and NLS type.
{\thm\label{thm_reg_DX} Let $n\geqslant2$ be an integer, and $\mathbf{e}=(\epsilon_1, \epsilon_2, \cdots, \epsilon_n)$ be a vector in the vector space $\bF^n_2$ and $\mathbf{a}=(\alpha_1, \alpha_2, \cdots, \alpha_n)$ be a vector in the vector space $\bF^n_4$. Let $G_\mathbf{e}$ and $Q_\mathbf{a}$ be the abelian group and quadratic form defined in Theorem~\ref{thm_DX}, and $\cS^{(4)}_{\mathbf{e},\mathbf{a}}$ be the amorphic abelian Cayley scheme in Theorem~\ref{thm_DX}. For each pair of integers $k$ and $l$ such that $0\leqslant l\leqslant w(\mathbf{e})$ and $1\leqslant n-l\leqslant k\leqslant n$, the abelian Cayley scheme admits a regular automorphism group $G$ of nilpotency class $2$ and of exponent $4$ such that $[G,G]\cong\bZ_2^{2k}$ or $\bZ_2^{2k-1}$, $Z(G)\cong\bZ_2^{4n-2k-4l}\oplus\bZ_4^{2l}$, or $\bZ_2^{4n-2k-4l+1}\oplus\bZ_4^{2l-1}$ or $\bZ_2^{4n-2k-4l-1}\oplus\bZ_4^{2l}$, and $\Phi(G)\cong\bZ_2^{2l+2w(\mathbf{e})-1}$ or $\bZ_2^{2l+2w(\mathbf{e})}$. }
\proof Let $\mathbf{v}=(\nu_1,\nu_2,\cdots,\nu_n)\in\bF^n_2$ be a vector such that $w(\mathbf{v}*\mathbf{e}+\mathbf{e})=l$  and $w(\mathbf{v})=k$, we define an automorphism $\tau_\mathbf{v}=(\tau_{\nu_1},\tau_{\nu_2},\cdots,\tau_{\nu_n})\in\Aut(G_\mathbf{e})$ of $G_\mathbf{e}$, which is also an isometry of $Q_\mathbf{a}$ of order $2$. Since $${\rm Im}_{G_\mathbf{e}}(1+\tau_\mathbf{v})={\rm Im}_{G_{\epsilon_1}}(1+\tau_{\nu_1})\oplus{\rm Im}_{G_{\epsilon_2}}(1+\tau_{\nu_2})\oplus\cdots\oplus{\rm Im}_{G_{\epsilon_n}}(1+\tau_{\nu_n}),$$ by Lemma~\ref{lem_index2} and Eqn.~(\ref{o2}), fix a vector $(b_1,b_2,\cdots,b_n)\ne\mathbf{0}$ in $\bF^n_4$, the subgroup $$K:=\{((x_1,y_1),(x_2,y_2),\cdots,(x_n,y_n))\in G_\mathbf{e}: \Tr(b_1x_1+b_2x_2+\cdots+b_nx_n)=0\}$$ of $G_\mathbf{e}$ is a $\tau_\mathbf{v}$-invariant subgroup of index $2$. Let $h=((x_1,y_1),(x_2,y_2),\cdots,(x_n,y_n))$ be an element in $G_\mathbf{e}$ such that $\Tr(b_1x_1+b_2x_2+\cdots+b_nx_n)\ne 0$. By Theorem~\ref{thm_Order2} the group $G:=G_{K,\tau_\mathbf{v},h}$ defined in Theorem~\ref{thm_GKt} is a regular automorphism group of the amorphic abelian group scheme $\cS^{(4)}_{\mathbf{e},\mathbf{a}}$. By Theorem~\ref{thm_GKt} and Eqn. (\ref{o2}), the commutator subgroup
\begin{align*} [G,G]&\cong{\rm Im}_K(1-\tau_\mathbf{v})\\
                    &=\{((0,\nu_1x_1),(0,\nu_2x_2),\cdots,(0,\nu_nx_n))\in G_\mathbf{e}:\Tr(b_1x_1+b_2x_2+\cdots+b_nx_n)=0\}\\
                    &\cong\begin{cases}\bZ_2^{2w(\mathbf{v})-1} &\text{ if }\mathbf{b}*(\mathbf{v}+\mathbf{1})=\mathbf{0},\\ \bZ_2^{2w(\mathbf{v})} &\text{ if }\mathbf{b}*(\mathbf{v}+\mathbf{1})\ne\mathbf{0}.\end{cases}
\end{align*}
By Eqn. (\ref{eqn_center}), the center
\begin{align*} Z(G)&\cong{\rm Ker}_K(1-\tau_\mathbf{v})\\
                   &=\{(((\nu_1+1)x_1,y_1),((\nu_2+1)x_2,y_2),\cdots,((\nu_n+1)x_n,y_n))\in G_\mathbf{e}:\Tr(b_1(\nu_1+1)x_1\\
                   &\hspace{80mm}+b_2(\nu_2+1)x_2+\cdots+b_n(\nu_n+1)x_n)=0\}\\
                   &\cong\begin{cases}\bZ_2^{4n-2w(\mathbf{v})-4w(\mathbf{v}*\mathbf{e}+\mathbf{e})}\oplus\bZ_4^{2w(\mathbf{v}*\mathbf{e}+\mathbf{e})} &\text{ if }\mathbf{b}*(\mathbf{v}+\mathbf{1})=\mathbf{0},\\ \bZ_2^{4n-2w(\mathbf{v})-4w(\mathbf{v}*\mathbf{e}+\mathbf{e})+1}\oplus\bZ_4^{2w(\mathbf{v}*\mathbf{e}+\mathbf{e})-1}  &\text{ if }\mathbf{b}*(\mathbf{v}+\mathbf{1})\ne\mathbf{0}, \mathbf{b}*(\mathbf{v}+\mathbf{1})*(\mathbf{e}+\mathbf{1})=\mathbf{0},\\
                   \bZ_2^{4n-2w(\mathbf{v})-4w(\mathbf{v}*\mathbf{e}+\mathbf{e})-1}\oplus\bZ_4^{2w(\mathbf{v}*\mathbf{e}+\mathbf{e})}&\text{ if }\mathbf{b}*(\mathbf{v}+\mathbf{1})*(\mathbf{e}+\mathbf{1})\ne\mathbf{0}.\\
                   \end{cases}
\end{align*}
By Eqn. (\ref{eqn_Frat}), the Frattini subgroup
\begin{align*} \Phi(G)&\cong\Phi(K)+{\rm Im}_{K}(1+\tau_\mathbf{v})+\<h_2\>=\Phi(K)+{\rm Im}_{K}(1+\tau_\mathbf{v})+\<h+\tau_\mathbf{v}(h)\>\\
                      &=\Phi(K)+{\rm Im}_{G_\mathbf{e}}(1+\tau_\mathbf{v})\\
                      &=\{((0,(\nu_1+\epsilon_1)x_1+\epsilon_1x'_1),(0,(\nu_2+\epsilon_2)x_2+\epsilon_2x'_2),\cdots,(0,(\nu_n+\epsilon_n)x_n+\epsilon_nx'_n))\in G_\mathbf{e}:\\ &\hspace{100mm}\Tr(b_1x'_1+b_2x'_2+\cdots+b_nx'_n)=0\}\\
                      &\cong\begin{cases}\bZ_2^{2w(\mathbf{v}*\mathbf{e}+\mathbf{v})+2w(\mathbf{e})-1}&\text{ if }\mathbf{b}*(\mathbf{v}*\mathbf{e}+\mathbf{1})=\mathbf{0},\\
                      \bZ_2^{2w(\mathbf{v}*\mathbf{e}+\mathbf{v})+2w(\mathbf{e})}&\text{ if }\mathbf{b}*(\mathbf{v}*\mathbf{e}+\mathbf{1})\ne\mathbf{0}.
                      \end{cases}
\end{align*}
Since ${\rm Im}_K(1-\tau_\mathbf{v})\ne\{0\}$ and ${\rm Im}_{\Phi(K)}(1-\tau_\mathbf{v})=\{0\}$, by Theorem~\ref{thm_Order2}, the nilpotency class of $G$ is $2$. It is easy to check that the exponent of $G$ is $4$. \qed

{\thm Let $n\geqslant2$ be an integer, and both $\mathbf{e}=(\epsilon_1, \epsilon_2, \cdots, \epsilon_n)$ and $\mathbf{a}=(\alpha_1, \alpha_2, \cdots, \alpha_n)$ be vectors in $\bF^n_2$. Let $G_\mathbf{e}$ and $Q_\mathbf{a}$ be the abelian group and quadratic form defined in Theorem~\ref{thm_DX}, and $\cS^{(3)}_{\mathbf{e},\mathbf{a}}$ be the amorphic abelian Cayley scheme in Theorem~\ref{thm_DX}. The abelian group scheme admits a regular automorphism group $G$ of nilpotency class $2$ or $3$ and of exponent $4$ or $8$ such that $[G,G]\cong\bZ_2^{2(n-w(\mathbf{e}))+1}\oplus\bZ_4^{w(\mathbf{e})-1}$ or $\bZ_2^{2(n-w(\mathbf{e}))-1}\oplus\bZ_4^{w(\mathbf{e})}$, $Z(G)\cong\bZ_2^{2(n-w(\mathbf{e}))}\oplus\bZ_4^{w(\mathbf{e})}$, and $\Phi(G)\cong\bZ_2^{2n-w(\mathbf{e})-1}\oplus\bZ_4^{w(\mathbf{e})}$ or $\bZ_2^{2n-w(\mathbf{e})}\oplus\bZ_4^{w(\mathbf{e})}$.}
\proof We define an automorphism $\rho_{\mathbf{a}}=(\rho_{\alpha_1},\rho_{\alpha_2},\cdots,\rho_{\alpha_n})$ of the abelian group $G_\mathbf{e}$, which is also a generalized isometry of $Q_\mathbf{a}$, of order $2$. Fix a non-zero vector $\mathbf{b}=(b_1,b_2,\cdots,b_n)\in\bF^n_2$. From
$${\rm Im}_{G_\mathbf{e}}(1+\rho_\mathbf{a})={\rm Im}_{G_{\epsilon_1}}(1+\rho_{\alpha_1})\oplus{\rm Im}_{G_{\epsilon_2}}(1+\rho_{\alpha_2})\oplus\cdots\oplus{\rm Im}_{G_{\epsilon_n}}(1+\rho_{\alpha_n})$$ and Eqn.~(\ref{image_1+rho}), by Lemma~\ref{lem_index2}, the subgroup $$K:=\{((x_1,y_1),(x_2,y_2),\cdots,(x_n,y_n))\in G_\mathbf{e} : \Tr(b_1x_1+b_2x_2+\cdots+b_nx_n)=0\}$$ is a $\rho_\mathbf{a}$-invariant subgroup of index $2$ in $G_\mathbf{e}$. Let $h$ be any element in $G_\mathbf{e}\setminus H$. The group $G:=G_{K,\rho_\mathbf{a},h}$ defined in Theorem~\ref{thm_GKt} is a regular automorphism group of the amorphic abelian Cayley scheme $\cS^{(3)}_{\mathbf{e},\mathbf{a}}$. By Theorem~\ref{thm_GKt} and Eqn. (\ref{image_1-rho}), the commutator subgroup
\begin{align*} [G,G]&\cong{\rm Im}_K(1-\rho_\mathbf{a})\\
                    &=\{((x_1+x^2_1,y_1+y^2_1+(\alpha_1+\epsilon_1)x^2_1+\epsilon_1x^3_1),(x_2+x^2_2,y_2+y^2_2+(\alpha_2+\epsilon_2)x^2_2+\epsilon_2x^3_2),\cdots,\\
                     &\hspace{10mm}(x_n+x^2_n,y_n+y^2_n+(\alpha_n+\epsilon_n)x^2_n+\epsilon_nx^3_n))\in G_\mathbf{e}: \Tr(b_1x_1+b_2x_2+\cdots+b_nx_n)=0\}\\
                    &\cong\begin{cases}\bZ_2^{2(n-w(\mathbf{e}))+1}\oplus\bZ_4^{w(\mathbf{e})-1}&\text{ if }\mathbf{b}*\mathbf{e}=\mathbf{b}\\ \bZ_2^{2(n-w(\mathbf{e}))-1}\oplus\bZ_4^{w(\mathbf{e})}&\text{ if }\mathbf{b}*\mathbf{e}\ne\mathbf{b}. \end{cases}
\end{align*}
By Eqn. (\ref{image_1-rho}), the center
\begin{align*} Z(G)&\cong{\rm Ker}_K(1-\rho_\mathbf{a})\\
                   &=\{((x_1,y_1),(x_2,y_2),\cdots,(x_n,y_n))\in G_\mathbf{e}: \Tr(b_1x_1+b_2x_2+\cdots+b_nx_n)=0, x_i+x_i^2=0\\
                   &\hspace{80mm}\text{ and }y_i+y_i^2=\alpha_ix_i^2\text{ for }i=1,2,\cdots,n\}\\
                   &\cong\bZ_2^{2(n-w(\mathbf{e}))}\oplus\bZ_4^{w(\mathbf{e})}.
\end{align*}
By Eqn. (\ref{eqn_Frat}), the Frattini subgroup
\begin{align*} \Phi(G)&\cong\Phi(K)+{\rm Im}_{K}(1+\rho_\mathbf{a})+\<h_2\>=\Phi(K)+{\rm Im}_{K}(1+\rho_\mathbf{a})+\<h+\rho_\mathbf{a}(h)\>\\
                      &=\Phi(K)+{\rm Im}_{G_\mathbf{e}}(1+\rho_\mathbf{a})\\
                      &=\{((x_1+x_1^2,y_1+y_1^2+\alpha_1x^2_1+\epsilon_1x^3_1+\epsilon_1x'_1),(x_2+x_2^2,y_2+y_2^2+\alpha_2x^2_2+\epsilon_2x^3_2+\epsilon_2x'_2),\cdots,\\
                       &\hspace{10mm}(x_n+x_n^2,y_n+y_n^2+\alpha_nx^2_n+\epsilon_nx^3_n+\epsilon_nx'_n))\in G_\mathbf{e}: \Tr(b_1x'_1+b_2x'_2+\cdots+b_nx'_n)=0\}\\
                      &\cong\begin{cases}\bZ_2^{2n-w(\mathbf{e})-1}\oplus\bZ_4^{w(\mathbf{e})}&\text{ if }\mathbf{b}*\mathbf{e}=\mathbf{b}\\
                      \bZ_2^{2n-w(\mathbf{e})}\oplus\bZ_4^{w(\mathbf{e})}&\text{ if }\mathbf{b}*\mathbf{e}\ne\mathbf{b}.
                      \end{cases}
\end{align*}
From the calculation of $[G,G]$, we know that $G$ is non-abelian. When $\mathbf{e}=\mathbf{0}$ or $w(\mathbf{e})=w(\mathbf{b})=w(\mathbf{e}*\mathbf{b})=1$, the restriction $\rho_\mathbf{a}|_{\Phi(K)}=1_{\Phi(K)}$ and $G$ is of nilpotency class $2$, and otherwise $G$ is of nilpotency class $3$. If $\mathbf{e}=\mathbf{0}$, the exponent of $G$ is $4$, and if $\mathbf{e}\ne\mathbf{0}$, the exponent of $G$ is $8$. \qed

{\thm Let $n\geqslant2$ be an integer, and $\mathbf{e}=(\epsilon_1, \epsilon_2, \cdots, \epsilon_n)$ be a vector in $\bF^n_2$ and $\mathbf{a}=(\alpha_1, \alpha_2, \cdots, \alpha_n)$ be a vector in $\bF^n_4\setminus\bF^n_2$. Let $G_\mathbf{e}$ and $Q_\mathbf{a}$ be the abelian group and quadratic form defined in Theorem~\ref{thm_DX}, and $\cS^{(3)}_{\mathbf{e},\mathbf{a}}$ be the amorphic abelian group scheme in Theorem~\ref{thm_DX}. Then the scheme admits a non-abelian regular automorphism group $G$ whose nilpotency class is $2$ or $4$ and whose exponent is $4$ or $8$ such that $|[G,G]|=2^{2(n-1)+w(\Tr(\mathbf{a}))}$ or $2^{2n+w(\Tr(\mathbf{a}))-1}$, $|Z(G)|=2^{2n-w(\Tr(\mathbf{a}))}$ or $2^{2n-w(\Tr(\mathbf{a}))-1}$, and $|\Phi(G)|=2^{2n-1+w(\Tr(\mathbf{a}))+w((\Tr(\mathbf{a})+\mathbf{1})*\mathbf{e})}$. }
\proof We define an automorphism $\rho_{\mathbf{a}}=(\rho_{\alpha_1},\rho_{\alpha_2},\cdots,\rho_{\alpha_n})$ of the abelian group $G_\mathbf{e}$, which is also a generalized isometry of $Q_\mathbf{a}$ of order $4$. Fix a non-zero vector $\mathbf{b}=(b_1,b_2,\cdots,b_n)\in\bF^n_2$. From
$${\rm Im}_{G_\mathbf{e}}(1+\rho_\mathbf{a})={\rm Im}_{G_{\epsilon_1}}(1+\rho_{\alpha_1})\oplus{\rm Im}_{G_{\epsilon_2}}(1+\rho_{\alpha_2})\oplus\cdots\oplus{\rm Im}_{G_{\epsilon_n}}(1+\rho_{\alpha_n})$$ and Eqn.~(\ref{image_1+rho}), by Lemma~\ref{lem_index2}, the subgroup $$H:=\{((x_1,y_1),(x_2,y_2),\cdots,(x_n,y_n))\in G_\mathbf{e} : \Tr(b_1x_1+b_2x_2+\cdots+b_nx_n)=0\}$$ is a $\rho_\mathbf{a}$-invariant subgroup of index $2$ in $G_\mathbf{e}$ and $$K:=\{((x_1,y_1),(x_2,y_2),\cdots,(x_n,y_n))\in G_\mathbf{e} : b_1x_1+b_2x_2+\cdots+b_nx_n=0\}$$ is a $\rho_\mathbf{a}$-invariant subgroup of index $2$ in $H$. Let $h$ be any element in $G_\mathbf{e}\setminus H$. The group $G:=G_{K,\rho_\mathbf{a},h}$ defined in Theorem~\ref{thm_GKt} is a regular subgroup of both graphs $\Gamma^{(0)}_{\mathbf{e},\mathbf{a}}$ and $\Gamma^{(1)}_{\mathbf{e},\mathbf{a}}$, corresponding to the partial difference sets $Q_{\mathbf{a}}^{-1}(0)\setminus \{0\}$ and  $Q_{\mathbf{a}}^{-1}(1)$, respectively.
By Theorem~\ref{thm_GKt} and Eqn. (\ref{image_1-rho}), the commutator subgroup
\begin{align*} [G,G]&\cong{\rm Im}_K(1-\rho_\mathbf{a})\\
                    &=\{((x_1+x^2_1,y_1+y^2_1+(\alpha_1+\epsilon_1)x^2_1+\epsilon_1x^3_1),(x_2+x^2_2,y_2+y^2_2+(\alpha_2+\epsilon_2)x^2_2+\epsilon_2x^3_2),\cdots,\\
                     &\hspace{20mm}(x_n+x^2_n,y_n+y^2_n+(\alpha_n+\epsilon_n)x^2_n+\epsilon_nx^3_n))\in G_\mathbf{e}: b_1x_1+b_2x_2+\cdots+b_nx_n=0\}\\
                    &\cong\begin{cases}\bZ_2^{2(n-w(\mathbf{e}))+w(\Tr(\mathbf{a}))}\oplus\bZ_4^{w(\mathbf{e})-1} &\text{ if }\mathbf{b}*\Tr(\mathbf{a})=\mathbf{b}\text{ and }\mathbf{b}*\mathbf{e}=\mathbf{b}\\
                    \bZ_2^{2(n-w(\mathbf{e})-1)+w(\Tr(\mathbf{a}))}\oplus\bZ_4^{w(\mathbf{e})} &\text{ if }\mathbf{b}*\Tr(\mathbf{a})=\mathbf{b}\text{ and }\mathbf{b}*\mathbf{e}\ne\mathbf{b}\\
                    \bZ_2^{2(n-w(\mathbf{e}))+w(\Tr(\mathbf{a}))+1}\oplus\bZ_4^{w(\mathbf{e})-1} &\text{ if }\mathbf{b}*\Tr(\mathbf{a})\ne\mathbf{b}\text{ and }\mathbf{b}*\mathbf{e}=\mathbf{b}\\\bZ_2^{2(n-w(\mathbf{e}))+w(\Tr(\mathbf{a}))-1}\oplus\bZ_4^{w(\mathbf{e})}&\text{ if }\mathbf{b}*\Tr(\mathbf{a})\ne\mathbf{b}\text{ and }\mathbf{b}*\mathbf{e}\ne\mathbf{b}. \end{cases}
\end{align*}
Hence $|[G,G]|=2^{2(n-1)+w(\Tr(\mathbf{a}))}$ or $2^{2n+w(\Tr(\mathbf{a}))-1}$. By Theorem~\ref{thm_GKt} and Eqn. (\ref{kernel_1-rho}), if $$w(\Tr(\mathbf{a}))=w(\mathbf{b})=w(\Tr(\mathbf{a})*\mathbf{b})=1,$$ then $o(\rho_\mathbf{a}|_K)=2$ and the center
$Z(G)\cong{\rm Ker}_K(1-\rho_\mathbf{a})\times\bZ_2$ when $h_4\in\Phi({\rm Ker}_K(1-\rho_\mathbf{a}))$, and the center $Z(G)\cong[{\rm Ker}_K(1-\rho_\mathbf{a})/\<h_4\>]\times\bZ_4$ when $h_4\not\in\Phi({\rm Ker}_K(1-\rho_\mathbf{a}))$. Otherwise $o(\rho_\mathbf{a}|_K)=o(\rho_\mathbf{a})=4$ and the center $Z(G)\cong{\rm Ker}_K(1-\rho_\mathbf{a})$, where
\begin{align*} &{\rm Ker}_K(1-\rho_\mathbf{a})\\
                   =&\{((x_1,y_1),(x_2,y_2),\cdots,(x_n,y_n))\in G_\mathbf{e}:b_1x_1+b_2x_2+\cdots+b_nx_n=0, x_i+x_i^2=0\text{ and }\\
                   &\hspace{90mm} y_i+y_i^2=\alpha_ix_i^2\text{ for }i=1,2,\cdots,n\}\\
                   \cong&\begin{cases}\bZ_2^{w(\Tr(\mathbf{a}))+2w((\Tr(\mathbf{a})+\mathbf{1})*(\mathbf{e}+\mathbf{1}))}\oplus\bZ_4^{w((\Tr(\mathbf{a})+\mathbf{1})*\mathbf{e})}&\text{ if }\mathbf{b}*\Tr(\mathbf{a})=\mathbf{b}\\
                   \bZ_2^{w(\Tr(\mathbf{a}))+2w((\Tr(\mathbf{a})+\mathbf{1})*(\mathbf{e}+\mathbf{1}))+1}\oplus\bZ_4^{w((\Tr(\mathbf{a})+\mathbf{1})*\mathbf{e})-1} &\text{ if }\mathbf{b}*\Tr(\mathbf{a})\ne\mathbf{b}\text{ and }\mathbf{b}*\mathbf{e}=\mathbf{b}\\ \bZ_2^{w(\Tr(\mathbf{a}))+2w((\Tr(\mathbf{a})+\mathbf{1})*(\mathbf{e}+\mathbf{1}))-1}\oplus\bZ_4^{w((\Tr(\mathbf{a})+\mathbf{1})*\mathbf{e})} &\text{ if }\mathbf{b}*\Tr(\mathbf{a})\ne\mathbf{b}\text{ and }\mathbf{b}*\mathbf{e}\ne\mathbf{b}.
                   \end{cases}
\end{align*}
Hence $|Z(G)|=2^{2n-w(\Tr(\mathbf{a}))}$ or $2^{2n-w(\Tr(\mathbf{a}))-1}$.
By Theorem~\ref{thm_GKt}, the Frattini subgroup $\Phi(G)\cong\<\Phi(K)+{\rm Im}_K(1+\rho_\mathbf{a}),h_2\rho_\mathbf{a}^2\>$. Since
\begin{align*} &\Phi(K)+{\rm Im}_K(1+\rho_\mathbf{a})\\
              =&\{((x_1+x^2_1,y_1+y^2_1+\alpha_1x_1^2+\epsilon_1x^3_1+\epsilon_1x'_1),(x_2+x^2_2,y_2+y^2_2+\alpha_2x_2^2+\epsilon_2x^3_2+\epsilon_2x'_2),\cdots,\\
              &\hspace{20mm}(x_n+x^2_n,y_n+y^2_n+\alpha_nx^2_n+\epsilon_nx^3_n+\epsilon_nx'_n))\in G_\mathbf{e}: b_1x_1+b_2x_2+\cdots+b_nx_n=0,\\
              &\hspace{100mm}b_1x'_1+b_2x'_2+\cdots+b_nx'_n=0\}\\
              &\cong\begin{cases}\bZ_2^{2(n-w(\mathbf{e}))+w(\Tr(\mathbf{a}))+w((\Tr(\mathbf{a})+\mathbf{1})*\mathbf{e})}\oplus\bZ_4^{w(\mathbf{e})-1} &\text{ if }\mathbf{b}*\mathbf{e}=\mathbf{b}\\
                    \bZ_2^{2(n-w(\mathbf{e})-1)+w(\Tr(\mathbf{a}))+w((\Tr(\mathbf{a})+\mathbf{1})*\mathbf{e})}\oplus\bZ_4^{w(\mathbf{e})} &\text{ if }\mathbf{b}*\mathbf{e}\ne\mathbf{b}
                     \end{cases}
\end{align*} and $(h_2\rho_\mathbf{a}^2)^2=h_4=(1+\rho_\mathbf{a}+\rho_\mathbf{a}^2+\rho_\mathbf{a}^3)(h)=(1+\rho_\mathbf{a})(1+\rho_\mathbf{a}^2)(h)\in{\rm Im}_K(1+\rho_\mathbf{a})$ as $(1+\rho_\mathbf{a}^2)(h)\in K$, the order
$|\Phi(G)|=2^{2n-1+w(\Tr(\mathbf{a}))+w((\Tr(\mathbf{a})+\mathbf{1})*\mathbf{e})}$. From $|[G,G]|\ne0$, the nilpotency class of $G$ is greater than $1$. Since $(1-\rho_\mathbf{a})^4=2(1+\rho_\mathbf{a}^2)=0$ in ${\rm End}(G_\mathbf{e})$, the nilpotency class of $G$ is no greater than $4$. Also from $(1-\rho_\mathbf{a})^2=1+\rho_\mathbf{a}^2-2\rho_\mathbf{a}$ and $(1-\rho_\mathbf{a})^3=(1+\rho_\mathbf{a})(1-\rho_\mathbf{a}^2)$, the homomorphism $(1-\rho_\mathbf{a})^3|_K=0$ if and only if $w(\Tr(\mathbf{a}))=w(\mathbf{b})=w(\Tr(\mathbf{a})*\mathbf{b})=1$ which implies $(1-\rho_\mathbf{a})^2|_K=0$. Hence $G$ is of nilpotency class $2$ when $o(\rho_\mathbf{a}|_K)=2$ and of nilpotency class $4$ when $o(\rho_\mathbf{a}|_K)=4$. If $\Tr(\mathbf{a})+\mathbf{e}=\mathbf{0}$, then the exponent of $G$ is $4$. Otherwise the exponent of $G$ is $8$.
\qed

{\thm Let $n\geqslant0$ be an integer, $\mathbf{e}=(\epsilon_1,\epsilon_2,\cdots,\epsilon_n)$ be a vector in the vector space $\bF_2^{n}$ and $\mathbf{a}=(\alpha_1, \alpha_2,\cdots, \alpha_n)$ be a vector in the vector space $\bF_4^{n}$. Let $\mathbf{e'}=(\epsilon,\epsilon,\epsilon,\epsilon)$ for some element $\epsilon\in\bF_2$ and $\mathbf{a'}=(\alpha,\alpha,\alpha,\alpha)$ for some element $\alpha\in\bF_4$. Let $G_\mathbf{e'}$, $G_\mathbf{e}$ and $Q_\mathbf{a'}$, $Q_\mathbf{a}$ be the abelian group and quadratic form defined in Theorem~\ref{thm_DX}, and $G_{\mathbf{e'}\oplus\mathbf{e}}=G_\mathbf{e'}\oplus G_\mathbf{e}$ and $Q_{\mathbf{a'}\oplus\mathbf{a}}= Q_\mathbf{a'}\oplus Q_\mathbf{a}$.  Let $\cS^{(4)}_{\mathbf{e'}\oplus\mathbf{e},\mathbf{a'}\oplus\mathbf{a}}$ be the amorphic abelian Cayley scheme defined in Theorem~\ref{thm_DX}. For each integers $0\leqslant l\leqslant w(\mathbf{e})$ and $0\leqslant n-l\leqslant k\leqslant n$, the Cayley scheme $\cS^{(4)}_{\mathbf{e'}\oplus\mathbf{e},\mathbf{a'}\oplus\mathbf{a}}$ admits a regular automorphism group $G$ of nilpotency class $4$ or $6$ and exponent $8$ or $16$ such that $|[G,G]|=4^{5+k}$, $|Z(G)|=4^{n+l+2}$, and $|\Phi(G)|=2^{2(5+k+l)+1}$ or $2^{2(6+k+l)+1}$. }
\proof Let $\mathbf{v}=(\nu_1,\nu_2,\cdots,\nu_n)$ be a vector in $\bF_2^n$ and $\tau_\mathbf{v}\colon G_\mathbf{e}\rightarrow G_\mathbf{e}$ be the isometry of $Q_\mathbf{a}$ in $\Aut(G_\mathbf{e})$ defined in Proof of Theorem \ref{thm_reg_DX} with $w(\mathbf{v} * \mathbf{e} +\mathbf{ e}) = l$ and $w(\mathbf{v})=k$. We define an isometry $\pi\in\Aut(G_\mathbf{e'}\oplus G_\mathbf{e})$ of $Q_\mathbf{a'}\oplus Q_\mathbf{a}$ of order $4$ by
\begin{align*} \pi&\colon G_{\mathbf{e'}\oplus\mathbf{e}}\longrightarrow G_{\mathbf{e'}\oplus\mathbf{e}}\\
               \pi&(((x_1,y_1),(x_2,y_2),(x_3,y_3),(x_4,y_4)),z)=(((x_4,y_4),(x_1,y_1),(x_2,y_2),(x_3,y_3)),\tau_\mathbf{v}(z))\\ &\hspace{40mm}\text{ for all }((x_1,y_1),(x_2,y_2),(x_3,y_3),(x_4,y_4))\in G_\mathbf{e'}\text{ and }z\in G_\mathbf{e}.
\end{align*}
Let  $$K:=\{(((x_1,y_1),(x_2,y_2),(x_3,y_3),(x_4,y_4)),z)\in G_{\mathbf{e'}\oplus\mathbf{e'}}\colon \Tr(x_1)=\Tr(x_3)\text{ and } \Tr(x_2)=\Tr(x_4)\}$$ and $h=(((\omega,0),(0,0),(0,0),(0,0)),0)$. Then $K$ is a $\pi$-invariant subgroup of index $4$ of $G_{\mathbf{e'}\oplus\mathbf{e}}$ and the triple $K$, $\pi$ and $h$ satisfies conditions stated in Theorem \ref{thm_GKt}. Therefore $G:=G_{K,\pi,h}$ is a regular automorphism group of $\mathcal{S}^{(4)}_{\mathbf{e'}\oplus\mathbf{e}}$.
By Theorem~\ref{thm_GKt}, the commutator subgroup
\begin{align*} [G,G]&\cong{\rm Im}_K(1-\pi)=\{(((x_1+x_4,y_1+y_4+\epsilon x_4+\epsilon x_1^2x_4^2),(x_2+x_1,y_2+y_1+\epsilon x_1+\epsilon x_2^2x_1^2),\\
                    &\hspace{15mm}(x_3+x_2,y_3+y_2+\epsilon x_2+\epsilon x_3^2x_2^2),(x_4+x_3,y_4+y_3+\epsilon x_3+\epsilon x_4^2x_3^2)),z)\in G_{\mathbf{e'}\oplus\mathbf{e}}\colon \\
                    &\hspace{40mm}\Tr(x_1)=\Tr(x_3) \text{ and }\Tr(x_2)=\Tr(x_4)\text{ and }z\in{\rm Im}_{G_\mathbf{e}}(1-\tau_\mathbf{v})\}\\
                    &=\begin{cases}\bZ_2^{2(5+w(\mathbf{v}))}&\text{ if }\epsilon=0\\
                                   \bZ_2^{2(1+w(\mathbf{v}))}\oplus\bZ_4^4&\text{ if }\epsilon=1.
                    \end{cases}
\end{align*}
Since $o(\pi|_K)=4$, by Theorem~\ref{thm_GKt}, the center
\begin{align*} Z(G)\cong&{\rm Ker}_K(1-\pi)\\
                   =&\{(((x,y),(x,y),(x,y),(x,y)),z)\in G_{\mathbf{e'}\oplus\mathbf{e}}: x,y\in\bF_4,z\in\Ker_{G_\mathbf{e}}(1-\tau_\mathbf{v})\}\\
                   =&\begin{cases}\bZ_2^{2(w(\mathbf{v})+2w((\mathbf{v}+\mathbf{1})*(\mathbf{e}+\mathbf{1}))+2)}\oplus\bZ_4^{2w((\mathbf{v}+\mathbf{1})*\mathbf{e})}&\text{ if }\epsilon=0\\
                                  \bZ_2^{2(w(\mathbf{v})+2w((\mathbf{v}+\mathbf{1})*(\mathbf{e}+\mathbf{1})))}\oplus\bZ_4^{2(1+w((\mathbf{v}+\mathbf{1})*\mathbf{e}))}&\text{ if }\epsilon=1.
                   \end{cases}
\end{align*}
By Theorem~\ref{thm_GKt}, the Frattini subgroup $\Phi(G)\cong\<\Phi(K)+{\rm Im}_K(1+\pi),h_2\pi^2\>$ and
\begin{align*}      &\Phi(K)+{\rm Im}_K(1+\pi)\\
                   =&\{(((x_1+x_4,y_1+y_4+\epsilon x_1^2x_4^2+\epsilon x'_1),(x_2+x_1,y_2+y_1+\epsilon x_2^2x_1^2+\epsilon x'_2),\\
                    &(x_3+x_2,y_3+y_2+\epsilon x_3^2x_2^2+\epsilon x'_3),(x_4+x_3,y_4+y_3+\epsilon x_4^2x_3^2)+\epsilon x'_4),z)\in G_{\mathbf{e'}\oplus\mathbf{e}}\colon\\
                    &\Tr(x_1)=\Tr(x_3)\text{ and }\Tr(x_2)=\Tr(x_4), \Tr(x'_1)=\Tr(x'_3)\text{ and }\Tr(x'_2)=\Tr(x'_4)\\
                    &\text{ and }z\in{\rm Im}_{G_\mathbf{e}}(1+\tau_{\mathbf{v}})+\Phi(G_{\mathbf{e}})\}\\
                   =&\begin{cases}\bZ_2^{2(n+5-w((\mathbf{v}+\mathbf{1})*(\mathbf{e}+\mathbf{1}))}&\text{ if }\epsilon=0\\
                                  \bZ_2^{2(n+2-w((\mathbf{v}+\mathbf{1})*(\mathbf{e}+\mathbf{1}))}\oplus\bZ_4^{4}&\text{ if }\epsilon=1,
                   \end{cases}
\end{align*}
which follows that $|\Phi(G)|=2^{2(5+k+l)+1}$ or $2^{2(6+k+l)+1}$.
It is clear that the element $x=(((1,0),(0,0),(0,0),(0,0)),0)\in K$ and $$(1+\pi+\pi^2+\pi^3)(x)=(((1,0),(1,0),(1,0),(1,0)),0)\ne0.$$ When $\epsilon=0$, we find that the endomorphisms $$(1-\pi)^3=1+\pi+\pi^2+\pi^3\text{ and } (1-\pi)^4=2(1+\pi^2)=0$$ in ${\rm End}(G_{\mathbf{e'}\oplus\mathbf{e}})$ and order $o(h\pi)=8$, hence the nilpotency class of $G$ is $4$ and the exponent of $G$ is $8$. When $\epsilon=1$, we have $$(1-\pi)^5=2(1+\pi+\pi^2+\pi^3)\text{ and }(1-\pi)^5(x)=(((0,1),(0,1),(0,1),(0,1)),0)\ne0$$ and the order $o(h\pi)=16$, therefore the nilpotency class of $G$ is $6$ and the exponent of $G$ is $16$. \qed

\section{Concluding Remark}
In this paper, we constructed non-abelian regular automorphism groups of exponent $4$, $8$ and $16$ and of nilpotency class $2$, $3$, $4$, and $6$ for the Davis-Xiang graphs and RT2 graphs. Such a regular automorphism groups give rise to amorphic non-abelian Cayley association schemes and these groups contain partial difference set with the same parameter set as the RT2 as well as Davis-Xiang strongly regular graphs. It seems hard to consider the isomorphism groups amongst all the resulting regular subgroups, since testing the isomorphism between two $2$-groups is a hard question. Instead, we have only computed some group theoretical invariants.  Our results suggest that a good way to construct non-abelian groups that contain non-trivial partial difference sets and amorphic non-abelian Cayley association schemes is to consider the regular subgroups of the known strongly regular graphs and known amorphic association schemes.

\section*{Acknowledgement}
The work of the first two authors was supported by National Natural Science Foundation of China under   Grant No. 11771392.

\end{document}